\documentclass[12pt,reqno]{amsart}
\usepackage{enumerate, latexsym, amsmath, amsfonts, amssymb, amsthm, color}
\def\pmod #1{\ ({\rm{mod}}\ #1)}
\def\Z{\Bbb Z}
\def\N{\Bbb N}

\def\Q{\Bbb Q}

\def\l{\left}
\def\r{\right}
\def\bg{\bigg}
\def\({\bg(}
\def\){\bg)}
\def\t{\text}
\def\f{\frac}

\def\sign{{\rm sign}}

\def\ls{\leqslant}
\def\gs{\geqslant}

\def\bi{\binom}

\def\ve{\varepsilon}

\def\eq{\equiv}

\def\da{\delta}

\def\Proof{\noindent{\it Proof}}
\def\Ack{\medskip\noindent {\bf Acknowledgment}}
\theoremstyle{plain}
\newtheorem{theorem}{Theorem}

\newtheorem{lemma}{Lemma}

\theoremstyle{definition}

\theoremstyle{remark}
\newtheorem{remark}{Remark}

\allowdisplaybreaks

 \vspace{4mm}

\begin{document}

\hbox{Accepted by Electron. Res. Arch.}
\smallskip

\title
[{Proof of some conjectures involving quadratic residues}]
{Proof of some conjectures involving quadratic residues}

\author
[Fedor Petrov and Zhi-Wei Sun] {Fedor Petrov and Zhi-Wei Sun}

\address {(Fedor Petrov) St. Petersburg Department of Steklov Mathematical Institute of Russian Academy of Sciences, Fontanka 27, 191023, St. Petersburg, Russia}
\email{fedyapetrov@gmail.com}

\address{(Zhi-Wei Sun) Department of Mathematics, Nanjing
University, Nanjing 210093, People's Republic of China}
\email{zwsun@nju.edu.cn}

\keywords{Quadratic residues modulo primes, quadratic fields, roots of unity, permutations, triangular numbers.
\newline \indent The work is supported by the NSFC-RFBR Cooperation and Exchange Program (grants NSFC 11811530072 and
RFBR 18-51-53020-GFEN-a). The second author is also supported
by the Natural Science Foundation of China (grant 11971222).}

\begin{abstract}
We confirm several conjectures
of Sun involving quadratic residues modulo odd primes. For any
prime $p\equiv 1\pmod 4$ and integer $a\not\equiv0\pmod p$, we prove that
\begin{align*}&(-1)^{|\{1\ls k<\frac p4:\ (\frac kp)=-1\}|}\prod_{1\ls j<k\ls(p-1)/2}(e^{2\pi iaj^2/p}+e^{2\pi iak^2/p})
\\=&\begin{cases}1&\text{if}\ p\equiv1\pmod 8,\\\left(\frac ap\right)\varepsilon_p^{-(\frac ap)h(p)}&\text{if}\ p\equiv5\pmod8,\end{cases}
\end{align*}
and that \begin{equation*}
\begin{aligned}&\l|\l\{(j,k):\ 1\ls j<k\ls\f{p-1}2\ \&\ \{aj^2\}_p>\{ak^2\}_p\r\}\r|
\\&+\l|\l\{(j,k):\ 1\ls j<k\ls\f{p-1}2\ \&\ \{ak^2-aj^2\}_p>\f p2\r\}\r|
\\\eq&\l|\l\{1\ls k<\f p4:\ \l(\f kp\r)=\l(\f ap\r)\r\}\r|\pmod2.
\end{aligned}
\end{equation*}
where $(\frac{a}p)$ is the Legendre symbol, $\varepsilon_p$ and $h(p)$ are the fundamental unit and the class number of the real quadratic field $\mathbb Q(\sqrt p)$ respectively, and
$\{x\}_p$ is the least nonnegative residue of an integer $x$ modulo $p$.
Also, for any prime $p\equiv3\pmod4$ and $\da=1,2$, we determine
$$(-1)^{\left|\left\{(j,k): \ 1\ls j<k\ls(p-1)/2\ \text{and}\ \{\da T_j\}_p>\{\da T_k\}_p\right\}\right|},$$
where $T_m$ denotes the triangular number $m(m+1)/2$.
\end{abstract}
\maketitle

\section{Introduction}
\setcounter{lemma}{0}
\setcounter{theorem}{0}
\setcounter{corollary}{0}
\setcounter{remark}{0}
\setcounter{equation}{0}

Let $p$ be an odd prime. It is well known that the numbers
$$1^2,\ 2^2,\ \ldots,\ \l(\f{p-1}2\r)^2$$
are pairwise incongruent modulo $p$, and they give all the $(p-1)/2$ quadratic residues modulo $p$.
 Recently Z.-W. Sun \cite{S19} initiated the study of permutations related to quadratic residues modulo $p$ as well as evaluations of related products involving $p$th roots of unity; many of his results
 in \cite{S19} are related to the class number of the quadratic field $\Q(\sqrt{(-1)^{(p-1)/2}p})$.
 In this paper, we confirm some conjectures of Sun \cite{S19} in this new direction.

 Let $p>3$ be a prime and let $\zeta=e^{2\pi i/p}$. Let $a\in\Z$ with $p\nmid a$.
Recently, Z.-W. Sun \cite[Theorem 1.3(ii)]{S19} showed that
\begin{align*}&\prod_{1\ls j<k\ls(p-1)/2}(\zeta^{aj^2}-\zeta^{ak^2})
\\=&\begin{cases}\pm i^{(p-1)/4}p^{(p-3)/8}\ve_p^{(\f ap)h(p)/2}&\t{if}\ p\eq1\pmod4,
\\(-p)^{(p-3)/8}&\t{if}\ p\eq3\pmod8,
\\(-1)^{(p+1)/8+(h(-p)-1)/2}(\f ap)p^{(p-3)/8}i&\t{if}\ p\eq7\pmod8,
\end{cases}\end{align*}
where $(\f{a}p)$ is the Legendre symbol, $\ve_p$ with $p\eq1\pmod4$ is the fundamental unit of the real quadratic field $\Q(\sqrt{p})$, and $h(d)$ with $d\eq0,1\pmod4$ not a square is the class number of the quadratic field with discriminant $d$.
Sun \cite[Theorem 1.5]{S19} also proved that
\begin{equation}\label{plus}\begin{aligned}&(-1)^{a\f{p+1}2\lfloor\f{p-1}4\rfloor}2^{(p-1)(p-3)/8}\prod_{1\ls j<k\ls(p-1)/2}\cos\pi\f{a(k^2-j^2)}p
\\=&
\prod_{1\ls j<k\ls(p-1)/2}(\zeta^{aj^2}+\zeta^{ak^2})
=\begin{cases}1&\t{if}\ p\eq3\pmod4,
\\\pm\ve_p^{(\f ap)h(p)((\f2p)-1)/2}&\t{if}\ p\eq1\pmod4.\end{cases}
\end{aligned}\end{equation}

Our first theorem confirms \cite[Conjecture 6.7]{S19}.

\begin{theorem}
Let $p$ be a prime with $p\eq1\pmod4$, and let $\zeta=e^{2\pi i/p}$.
Let $a$ be an integer not divisible by $p$.

{\rm (i)} If $p\eq1\pmod 8$, then
\begin{equation}\label{main}\prod_{1\ls j<k\ls(p-1)/2}(\zeta^{aj^2}+\zeta^{ak^2})=(-1)^{|\{1\ls k<\f p4:\ (\f kp)=-1\}|}.
\end{equation}

{\rm (ii)} When $p\eq5\pmod 8$, we have
\begin{equation}\label{main5}(-1)^{|\{1\ls k<\f p4:\ (\f kp)=-1\}|}\prod_{1\ls j<k\ls(p-1)/2}(\zeta^{aj^2}+\zeta^{ak^2})
=\l(\f ap\r)\ve_p^{-(\f ap)h(p)}.\end{equation}
\end{theorem}
\begin{remark}\label{Rem1.1}
Let $p$ be a prime with $p\eq1\pmod4$. Then $(\f{p-1}2!)^2\eq-1\pmod p$ by Wilson's theorem. We may write $p=x^2+y^2$ with $x,y\in\Z$, $x\eq1\pmod4$ and $y\eq \f{p-1}2!x\pmod p$.
As $y^2\eq p-1\pmod 8$, we see that $y\eq(\f 2p)-1\pmod4$.
By a result of K. Burde \cite{B}, we have
$$\l|\l\{1\ls k< \f p4:\ \l(\f kp\r)=1\r\}\r|\eq0\pmod2\iff y\eq\l(\f 2p\r)-1\pmod 8.$$
Thus
\begin{equation}\label{yy}(-1)^{|\{1\ls k<\f p4:\ (\f kp)=-1\}|}=(-1)^{\f{p-1}4}(-1)^{\f14(y-(\f 2p)+1)}=(-1)^{\lfloor \f y4\rfloor}.
\end{equation}
\end{remark}
\medskip

Let $p$ be an odd prime. For each $a\in\Z$ we let $\{a\}_p$ denote the least nonnegative residue of $a$ modulo $p$. Define
$$s(p):=\l|\l\{(j,k):\ 1\ls j<k\ls\f{p-1}2\ \&\ \{j^2\}_p>\{k^2\}_p\r\}\r|$$
and
$$t(p):=\l|\l\{(j,k):\ 1\ls j<k\ls\f{p-1}2\ \&\ \{k^2-j^2\}_p>\f p2\r\}\r|$$
as in \cite{S19}, where $(j,k)$ is an ordered pair. Sun \cite[Theorem 1.4(i)]{S19} showed that
$$(-1)^{s(p)}=(-1)^{t(p)}=\begin{cases}1&\t{if}\ p\eq3\pmod 8,
\\(-1)^{(h(-p)+1)/2}&\t{if}\ p\eq7\pmod 8.
\end{cases}$$
He also conjectured
that (cf. \cite[Conjecture 6.1]{S19}) if $p\eq1\pmod4$ then
\begin{equation}\label{st}s(p)+t(p)\eq\l|\l\{1\ls k<\f p4:\ \l(\f kp\r)=1\r\}\r|\pmod2.
\end{equation}
Our second theorem in the case $a=1$ confirms this conjecture.

\begin{theorem}\label{Th1.2} Let $p$ be a prime with $p\eq1\pmod4$, and let $a\in\Z$ with $p\nmid a$. Then
\begin{equation}\label{parity}
\begin{aligned}&\l|\l\{(j,k):\ 1\ls j<k\ls\f{p-1}2\ \&\ \{aj^2\}_p>\{ak^2\}_p\r\}\r|
\\&+\l|\l\{(j,k):\ 1\ls j<k\ls\f{p-1}2\ \&\ \{ak^2-aj^2\}_p>\f p2\r\}\r|
\\\eq&\l|\l\{1\ls k<\f p4:\ \l(\f kp\r)=\l(\f ap\r)\r\}\r|\pmod2.
\end{aligned}
\end{equation}
\end{theorem}

Our third theorem was first conjectured by Sun (cf. \cite[Conjectures 6.3 and 6.4]{S19}).

\begin{theorem}\label{Th1.3} Let $p$ be a prime with $p\eq3\pmod4$.

{\rm (i)} We have
\begin{equation}\label{jk-sign}(-1)^{\l|\l\{(j,k): \ 1\ls j<k\ls(p-1)/2\ \t{and}\ \{j(j+1)\}_p>\{k(k+1)\}_p\r\}\r|}=(-1)^{\lfloor(p+1)/8\rfloor}.
\end{equation}

{\rm (ii)} Suppose $p>3$ and write $T_m=m(m+1)/2$ for $m\in\N$. Then
\begin{equation}\label{tri}\begin{aligned}(-1)^{|\{(j,k):\ 1\ls j<k\ls(p-1)/2\ \&\ \{T_j\}_p>\{T_k\}_p\}|}
=(-1)^{\f{h(-p)+1}2+|\{1\ls k\ls\lfloor\f{p+1}8\rfloor:\ (\f kp)=1\}|}.\end{aligned}
\end{equation}
\end{theorem}

We will prove Theorems 1.1-1.2 in Section 2.
Based on an auxiliary theorem given in Section 3,
we are going to prove Theorem 1.3 in Section 4.

\section{Proofs of Theorems 1.1-1.2}
\setcounter{lemma}{0}
\setcounter{theorem}{0}
\setcounter{corollary}{0}
\setcounter{remark}{0}
\setcounter{equation}{0}

In 2006, H. Pan \cite{P06} obtained the following lemma.

\begin{lemma}\label{Lem2.1} {\rm (H. Pan \cite{P06})} Let $n>1$ be an odd integer and let $c$ be any integer relatively prime to $n$. For each $j=1,\ldots,(n-1)/2$ let $\pi_c(j)$ be the unique $r\in\{1,\ldots,(n-1)/2\}$ with $cj$
congruent to $r$ or $-r$ modulo $n$. For the permutation $\pi_c$ on $\{1,\ldots,(n-1)/2\}$, its sign is given by
$$\sign(\pi_c)=\l(\f cn\r)^{(n+1)/2},$$
where $(\f cn)$ is the Jacobi symbol.
\end{lemma}

\medskip
\noindent{\it Proof of the First Part of Theorem 1.1}. As $p\eq1\pmod 8$, there is an integer $c$
with $c^2\eq2\pmod p$. For $j=1,\ldots,(p-1)/2$ let $\pi_c(j)$ be the unique $r\in\{1,\ldots,(p-1)/2\}$
with $cj$ congruent to $r$ or $-r$ modulo $p$. Then $\pi_c$ is a permutation on $\{1,\ldots,(p-1)/2\}$, and
\begin{align*}&
\prod_{1\ls j<k\ls(p-1)/2}\f{\zeta^{2aj^2}-\zeta^{2ak^2}}{\zeta^{aj^2}-\zeta^{ak^2}}
=\prod_{1\ls j<k\ls(p-1)/2}\f{\zeta^{a\pi_c(j)^2}-\zeta^{a\pi_c(k)^2}}{\zeta^{aj^2}-\zeta^{ak^2}}
\\=&(-1)^{|\{(j,k):\ 1\ls j<k\ls(p-1)/2\ \&\ \pi_c(j)>\pi_c(k)\}|}=\sign(\pi_c)=\l(\f cp\r)
\end{align*}
with the aid of Lemma 2.1. In view of K. S. Williams and J. D. Currie \cite[(1.4)]{WC}, we have
$$\l(\f cp\r)\eq c^{(p-1)/2}=(c^2)^{(p-1)/4}\eq 2^{(p-1)/4}\eq(-1)^{|\{0<k<\f p4:\ (\f kp)=-1\}|}\pmod p.$$
Therefore (\ref{main}) holds in the case $p\eq1\pmod 8$. \qed

\begin{remark}\label{Rem2.1} Our method to prove part (i) of Theorem 1.1 does not work for part (ii) of Theorem 1.1.
\end{remark}

\medskip
\noindent{\it Proof of the Second Part of Theorem 1.1}. We distinguish two cases.
\medskip

{\it Case} 1. $(\f ap)=1$.

In this case,
$$\l\{\{aj^2\}_p:\ 1\ls j\ls\f{p-1}2\r\}=\l\{\{k^2\}_p:\ 1\ls k\ls\f{p-1}2\r\}$$
So it suffices to show \eqref{main5} for $a=1$. In view of (1.1) with $a=1$,
we only need to prove that
\begin{equation}\label{2.1}(-1)^{|\{0<k<\f p4:\ (\f kp)=-1\}|}\prod_{1\ls j<k\ls(p-1)/2}(\zeta^{j^2}+\zeta^{k^2})>0.
\end{equation}

As $(\f{-1}p)=1$, for each $1\ls j\ls(p-1)/2$ there is a unique integer $j_*\in\{1,\ldots,(p-1)/2\}$ such that $p-j^2\eq j_*^2\pmod p$. As $(\f 2p)=-1$, we have $j\not=j_*$. For any distinct $j,k\in\{1,\ldots,(p-1)/2\}$, we have $\zeta^{j^2}+\zeta^{k^2}\not=0$ (since $\zeta^{2j^2}\not=\zeta^{2k^2}$) and
$$(\zeta^{j^2}+\zeta^{k^2})(\zeta^{j_*^2}+\zeta^{k_*^2})
=(\zeta^{j^2}+\zeta^{k^2})(\zeta^{-j^2}+\zeta^{-k^2})=|\zeta^{j^2}+\zeta^{k^2}|^2>0;$$
also,
$$\{j,k\}=\{j_*,k_*\}\iff j_*=k\ \t{and}\ k_*=j\iff j_*=k.$$
For $1\ls j\ls(p-1)/2$, clearly
$$\zeta^{j^2}+\zeta^{j_*^2}=\zeta^{j^2}+\zeta^{-j^2}=2\cos 2\pi\f{j^2}p=2\cos 2\pi\f{j_*^2}p$$
and hence
$$\zeta^{j^2}+\zeta^{j_*^2}>0\iff \cos 2\pi\f{j^2}p>0\iff \{j^2\}_p<\f p4\ \t{or}\ \{j_*^2\}_p<\f p4.$$
Thus the sign of the product
$$\prod_{1\ls j<k\ls(p-1)/2\atop p\mid j^2+k^2}(\zeta^{j^2}+\zeta^{k^2})
=(-1)^{(p-1)/4}\prod_{1\ls j<j_*\ls(p-1)/2}(-\zeta^{j^2}-\zeta^{j_*^2})$$
is
$$(-1)^{(p-1)/4-|\{1\ls k<\f p4:\ (\f kp)=1\}|}=(-1)^{|\{1\ls k<\f p4:\ (\f kp)=-1\}|}.$$
This proves (\ref{2.1}).

\medskip
{\it Case} 2. $(\f ap)=-1$.

By the discussion in Case 1, we have
\begin{equation}\label{2.2}(-1)^{|\{0<k<\f p4:\ (\f kp)=-1\}|}\prod_{1\ls j<k\ls(p-1)/2}(\zeta^{j^2}+\zeta^{k^2})=\ve_p^{-h(p)}.\end{equation}
Let $\varphi_a$ be the element of the Galois group $\text{Gal}(\Q(\zeta)/\Q)$ with $\varphi_a(\zeta)=\zeta^a$. Then
$$\varphi_a(\sqrt p)=\varphi_a\(\sum_{x=0}^{p-1}\zeta^{x^2}\)=\sum_{x=0}^{p-1}\zeta^{ax^2}=\l(\f ap\r)\sqrt{p}=-\sqrt p$$
by the evaluation of quadratic Gauss sums (cf. \cite[pp.\,70-75]{IR}). Hence
$$\varphi_a(\ve_p^{-h(p)})=\l(\f{N(\ve_p)}{\ve_p}\r)^{-h(p)}=-\ve_p^{h(p)}$$
where $N(\ve_p)$ is the norm of $\ve_p$ with respect to the field extension $\Q(\zeta_p)/\Q$, and we have used the known results $N(\ve_p)=-1$ and $2\nmid h(p)$ (cf. \cite[p.\,185 and p.\,187]{Co}).
Thus, by applying the automorphism $\varphi_a$ to the identity \eqref{2.2}, we get
\begin{align*}&(-1)^{|\{0<k<\f p4:\ (\f kp)=-1\}|}\prod_{1\ls j<k\ls(p-1)/2}(\zeta^{aj^2}+\zeta^{ak^2})
\\=&\varphi_a(\ve_p^{-h(p)})=-\ve_p^{h(p)}=\l(\f ap\r)\ve_p^{-(\f ap)h(p)}.
\end{align*}

In view of the above, we have proven Theorem 1.1(ii). \qed

\begin{lemma}\label{Lem3.1} Let $p$ be an odd prime, and let $a\in\Z$ with $p\nmid a$.
Then
\begin{equation}\label{3.1}\begin{aligned}&\l|\l\{(j,k):\ 1\ls j<k\ls\f{p-1}2\ \&\ \{aj^2\}_p>\{ak^2\}_p\r\}\r|
\\&+\l|\l\{(j,k):\ 1\ls j<k\ls\f{p-1}2\ \&\ \{ak^2-aj^2\}_p>\f p2\r\}\r|
\\\eq&\l|\l\{(j,k):\ 1\ls j<k\ls\f{p-1}2\ \&\ |\{aj^2\}_p-\{ak^2\}_p|>\f p2\r\}\r|\pmod2.
\end{aligned}\end{equation}
\end{lemma}
\Proof. This can be easily checked by distinguishing the cases $\{aj^2\}_p<\{ak^2\}_p$ and $\{aj^2\}_p>\{ak^2\}_p$ for $1\ls j<k\ls (p-1)/2$. \qed

\medskip
\noindent{\it Proof of Theorem \ref{Th1.2}}.
In view of Lemma \ref{Lem3.1}, it suffices to show that
\begin{equation}\label{2.4}\begin{aligned}&\l|\l\{(j,k):\ 1\ls j<k\ls \f{p-1}2\ \&\ |\{aj^2\}_p-\{ak^2\}_p|>\f p2\r\}\r|
\\&\qquad\eq\l|\l\{1\ls k<\f p4:\ \l(\f {ak}p\r)=1\r\}\r|\pmod2.
\end{aligned}\end{equation}

As $(\f{-1}p)=1$, for each $1\ls j\ls(p-1)/2$ there is a unique integer $j_*\in\{1,\ldots,(p-1)/2\}$ such that $-j^2\eq j_*^2\pmod p$ and hence $-aj^2\eq aj_*^2\pmod p$. Clearly,
$|\{1\ls j\ls(p-1)/2:\ j=j_*\}|\ls1$.
Note that
$$|\{aj_*^2\}_p-\{ak_*^2\}_p|=|(p-\{aj^2\}_p)-(p-\{ak^2\}_p)|=|\{aj^2\}_p-\{ak^2\}_p|.$$
If $j$ and $k$ are distinct elements of $\{1,\ldots,(p-1)/2\}$, then
$\{j,k\}=\{j_*,k_*\}$ if and only if
$j_*=k$ and $k_*=j$.
Thus
\begin{align*}&\l|\l\{(j,k):\ 1\ls j<k\ls \f{p-1}2\ \&\ |\{aj^2\}_p-\{ak^2\}_p|>\f p2\r\}\r|
\\\eq&\f12\l|\l\{1\ls j\ls\f{p-1}2:\ |\{aj^2\}_p-\{aj_*^2\}|=|2\{aj^2\}_p-p|>\f p2\r\}\r|
\\=&\f12\l|\l\{1\ls j\ls\f{p-1}2:\ \{aj^2\}_p<\f p4\ \t{or}\ \{aj^2\}_p>\f 34p\r\}\r|
\\=&\f12\l|\l\{1\ls j\ls\f{p-1}2:\ \{aj^2\}_p<\f p4\r\}\r|
+\f12\l|\l\{1\ls j\ls\f{p-1}2:\ \{aj_*^2\}_p<\f p4\r\}\r|
\\=&\l|\l\{1\ls j\ls\f{p-1}2:\ \{aj^2\}_p<\f p4\r\}\r|
=\l|\l\{1\ls k<\f p4:\ \l(\f{ak}p\r)=1\r\}\r|
\pmod2.
\end{align*}
This proves the desired \eqref{2.4}. \qed

\section{An Auxiliary Theorem}
\setcounter{lemma}{0}
\setcounter{theorem}{1}
\setcounter{corollary}{0}
\setcounter{remark}{0}
\setcounter{equation}{0}

We first need a result of Sun \cite{S19}.

\begin{lemma}\label{Lem4.1} Let $p=2n+1$ be a prime with $n$ odd, and let $a\in\Z$ with $p\nmid a$. Then
\begin{equation}\label{4.1}\begin{aligned}&(-1)^{|\{(j,k):\ 1\ls j<k\ls n\ \&\ \{aj^2\}_p>\{ak^2\}_p\}|}
\\=&\begin{cases}1&\t{if}\ p\eq3\pmod8,\\(-1)^{(h(-p)+1)/2}(\f ap)&\t{if}\ p\eq7\pmod8.
\end{cases}\end{aligned}\end{equation}
\end{lemma}
\proof By Sun \cite[Theorem 1.4(ii)]{S19},
\begin{align*}&\prod_{1\ls j<k\ls(p-1)/2}\l(\cot\pi\f{aj^2}p-\cot\pi\f{ak^2}p\r)
\\=&\begin{cases} (2^{p-1}/p)^{(p-3)/8}&\t{if}\ p\eq3\pmod8,
\\(-1)^{(h(-p)+1)/2}(\f ap)(2^{p-1}/p)^{(p-3)/8}&\t{if}\ p\eq7\pmod8.\end{cases}
\end{align*}
This implies \eqref{4.1} since for any $1\ls j<k\ls(p-1)/2$ we have
$$\cot\pi\f{aj^2}p<\cot\pi\f{ak^2}p\iff \{aj^2\}_p>\{ak^2\}_p.$$
We are done. \qed

\begin{theorem}\label{Th4.1} Let $p=2n+1$ be a prime with $n$ odd, and let $a,b\in\{1,\ldots,p-1\}$.
Then
\begin{equation}\label{4.2}
\begin{aligned}&(-1)^{|\{(s,t):\ 0\ls t<s\ls n\ \&\ \{as^2-b\}_p>\{at^2-b\}_p\}|-|\{0<r<b:\ (\f rp)=(\f ap)\}|}
\\=&\begin{cases}1&\t{if}\ p\eq3\pmod8,
\\(-1)^{(h(-p)-1)/2}(\f ap)&\t{if}\ p\eq7\pmod8.\end{cases}
\end{aligned}\end{equation}
\end{theorem}
\proof Let $0\ls t<s\ls n$. By comparing $\{as^2\}_p$ and $\{at^2\}_p$ with $b$, we verify case by case that $$[\{as^2-b\}_p>\{at^2-b\}_p]+[\{as^2\}_p>\{at^2\}_p]$$ is odd if and only if
$$\{as^2\}_p\gs b>\{at^2\}_p\ \ \t{or}\ \ \{at^2\}_p\gs b>\{as^2\}_p,$$
where for an assertion $A$ we define
$$[A]=\begin{cases}1&\t{if}\ A\ \t{holds},\\0&\t{otherwise}.\end{cases}$$
Note that
\begin{align*}&|\{(s,t):\ 0\ls t<s\ls n,\ \{as^2\}_p\gs b>\{at^2\}_p\ \t{or}\ \{as^2\}_p<b\ls\{at^2\}_p\}|
\\=&\l|\l\{(r_1,r_2):\ 0\ls r_1<b\ls r_2\ls p-1\ \&\ \l(\f{ar_1}p\r),\l(\f{ar_2}p\r)\not=-1\r\}\r|
\\=&\l|\l\{0\ls r<b:\ \l(\f {ar}p\r)\not=-1\r\}\r|
\times\l(\f{p+1}2-\l|\l\{0\ls r<b:\ \l(\f {ar}p\r)\not=-1\r\}\r|\r)
\\\eq&\l|\l\{0\ls r<b:\ \l(\f {ar}p\r)\not=-1\r\}\r|=1+\l|\l\{0<r<b:\ \l(\f rp\r)=\l(\f ap\r)\r\}\r|\pmod2
\end{align*}
and
\begin{align*}&|\{(s,t):\ 0\ls t<s\ls n\ \&\ \{as^2\}_p>\{at^2\}_p\}|
\\=&\bi{n+1}2-|\{(s,t):\ 0\ls t<s\ls n\ \&\ \{as^2\}_p<\{at^2\}_p\}|
\end{align*}
with $\bi{n+1}2\eq\f{p+1}4\pmod2$.
Combining the above with \eqref{4.1}, we finally obtain \eqref{4.2}. \qed

\section{Proof of Theorem 1.3}
\setcounter{lemma}{0}
\setcounter{theorem}{0}
\setcounter{corollary}{0}
\setcounter{remark}{0}
\setcounter{equation}{0}

\begin{lemma}\label{Lem5.1} Let $p$ be a prime with $p\eq3\pmod4$.

{\rm (i) (Dirichlet (cf. \cite[p.\,238]{IR}))} If $p>3$ then
$$\l(2-\l(\f2p\r)\r)h(-p)=\sum_{k=1}^{(p-1)/2}\l(\f kp\r).$$

{\rm (ii) (B. C. Berndt and S. Chowla \cite{BC})} If $p\eq3\pmod8$, then
$\sum_{0<k<p/4}(\f kp)=0.$
If $p\eq7\pmod8$, then
$\sum_{p/4<k<p/2}(\f kp)=0.$
\end{lemma}

\medskip
\noindent{\it Proof of Theorem \ref{Th1.3}}. We just prove the second part in details since the first part can be proved similarly.

Write $n=(p-1)/2$, and set
$$a=\f{p+1}2\ \ \t{and}\ \ b=\begin{cases}(5p+1)/8&\t{if}\ p\eq3\pmod8,
\\(p+1)/8&\t{if}\ p\eq7\pmod8.\end{cases}$$
For any $r\in\Z$, we have
$$T_{n-r}=\f{n(n+1)}2-(2n+1)\f r2+\f{r^2}2\eq ar^2-b\pmod p.$$
Thus
\begin{align*}&|\{(j,k):\ 0\ls j<k\ls n\ \&\ \{T_j\}_p>\{T_k\}_p\}|
\\=&|\{(t,s):\ 0\ls t<s\ls n\ \&\ \{T_{n-s}\}_p>\{T_{n-t}\}_p\}|
\\=&|\{(t,s):\ 0\ls t<s\ls n\ \&\ \{as^2-b\}_p>\{at^2-b\}_p\}|.
\end{align*}
Note that $(\f ap)=(\f 2p)$. Set
\begin{equation}\label{B}B:=\l|\l\{0<r<b:\ \l(\f rp\r)=\l(\f 2p\r)\r\}\r|.
\end{equation}
Applying Theorem \ref{Th4.1}, from the above we obtain
\begin{equation}\label{jkt}\begin{aligned}&|\{(j,k):\ 1\ls j<k\ls n\ \&\ \{T_j\}_p>\{T_k\}_p\}|
\\\eq&B+\begin{cases}0\pmod2&\t{if}\ p\eq3\pmod8,
\\(h(-p)-1)/2\pmod2&\t{if}\ p\eq7\pmod8.
\end{cases}\end{aligned}\end{equation}

When $p\eq7\pmod8$, we have
$$B+1=\l|\l\{1\ls k\ls\f{p+1}8:\ \l(\f kp\r)=1\r\}\r|$$
and hence \eqref{tri} follows from \eqref{jkt}.

Below we handle the case $p\eq3\pmod8$. Observe that
\begin{align*}B=&\sum_{k=1}^{(p-1)/2}\f{1-(\f kp)}2+\l|\l\{\f p2<k<\f{5p+1}8:\ \l(\f {2k-p}p\r)=1\r\}\r|
\\=&\f{p-1}4-\f12\sum_{k=1}^{(p-1)/2}\l(\f kp\r)+\l|\l\{0<r<\f{p+1}4:\ 2\nmid r\ \&\ \l(\f rp\r)=1\r\}\r|.
\end{align*}
Applying Lemma \ref{Lem5.1}, we obtain
\begin{align*}
B=&\f{p-1}4-\f{3h(-p)}2+\sum_{0<k<p/4}\f{1+(\f kp)}2
\\&-\l|\l\{0<r<\f{p+1}4:\ 2\mid r\ \&\ \l(\f rp\r)=1\r\}\r|
\\\eq&\f{h(-p)+1}2+\f{p-3}8-\l|\l\{0<k<\f{p+1}8:\ \l(\f {2k}p\r)=1\r\}\r|
\\=&\f{h(-p)+1}2+\l|\l\{0<k<\f{p+1}8:\ \l(\f {2k}p\r)=-1\r\}\r|
\\=&\f{h(-p)+1}2+\l|\l\{0<k\ls\l\lfloor\f{p+1}8\r\rfloor:\ \l(\f {k}p\r)=1\r\}\r|\pmod2.
\end{align*}
So, in this case, \eqref{tri} also follows from \eqref{jkt}. \qed
\medskip

\Ack. We would like to thank the referee for helpful comments.

\end{document}